\documentclass[12pt]{amsart}

\usepackage{epsfig}
\usepackage{subfig}
\usepackage{amscd}
\usepackage[mathscr]{eucal}
\usepackage{amssymb}
\usepackage{amsxtra}
\usepackage{amsmath}
\usepackage[all]{xy}
\usepackage{tikz}
\usepackage{pgf}
\usetikzlibrary{calc,backgrounds,arrows,matrix,shapes}
\usepackage{verbatim}
\usepackage[bookmarks]{hyperref}
\usepackage{mathtools}
\usepackage{tikz-cd}
\makeatletter

\pgfdeclareshape{genuspic}{
 \anchor{center}{\pgfpointorigin}
\backgroundpath{
     \pgfpathmoveto{\pgfqpoint{-1cm}{0cm}}
     \pgfpathcurveto %
            {\pgfpoint{-0.5cm}{-.5cm}}
        {\pgfpoint{0.5cm}{-.5cm}}
        {\pgfpoint{1cm}{0cm}}

     \pgfpathmoveto{\pgfqpoint{-0.75cm}{-0.15cm}}
     \pgfpathcurveto %
            {\pgfpoint{-0.25cm}{.25cm}}
        {\pgfpoint{.25cm}{.25cm}}
        {\pgfpoint{0.75cm}{-0.15cm}}
    }
    }

\pgfdeclareshape{strokegenuspic}{
 \anchor{center}{\pgfpointorigin}
\backgroundpath{
     \pgfpathmoveto{\pgfqpoint{-1cm}{0cm}}
     \pgfpathcurveto %
            {\pgfpoint{-0.5cm}{-.5cm}}
        {\pgfpoint{0.5cm}{-.5cm}}
        {\pgfpoint{1cm}{0cm}}

     \pgfpathmoveto{\pgfqpoint{-0.75cm}{-0.15cm}}
     \pgfpathcurveto %
            {\pgfpoint{-0.25cm}{.25cm}}
        {\pgfpoint{.25cm}{.25cm}}
        {\pgfpoint{0.75cm}{-0.15cm}}
        \pgfusepath{stroke}
    }
    }

\pgfdeclareshape{hackgenuspic}{
 \anchor{center}{\pgfpointorigin}
\backgroundpath{
         \pgfpathmoveto{\pgfqpoint{-0.78cm}{-.17cm}}
     \pgfpathcurveto %
            {\pgfpoint{-0.35cm}{-.44cm}}
        {\pgfpoint{0.35cm}{-.44cm}}
        {\pgfpoint{.78cm}{-0.17cm}} 
     \pgfpathmoveto{\pgfqpoint{-0.78cm}{-0.17cm}}
     \pgfpathcurveto %
            {\pgfpoint{-0.25cm}{.25cm}}
        {\pgfpoint{.25cm}{.25cm}}
        {\pgfpoint{0.78cm}{-0.17cm}}
        \pgfsetfillcolor{white}
        \pgfusepath{fill}

     \pgfpathmoveto{\pgfqpoint{-1cm}{0cm}}
     \pgfpathcurveto %
            {\pgfpoint{-0.5cm}{-.5cm}}
        {\pgfpoint{0.5cm}{-.5cm}}
        {\pgfpoint{1cm}{0cm}}

     \pgfpathmoveto{\pgfqpoint{-0.75cm}{-0.15cm}}
     \pgfpathcurveto %
            {\pgfpoint{-0.25cm}{.25cm}}
        {\pgfpoint{.25cm}{.25cm}}
        {\pgfpoint{0.75cm}{-0.15cm}}
              \pgfusepath{stroke}
    }
    }

    \pgfdeclareshape{fillhackgenuspic}{
 \anchor{center}{\pgfpointorigin}
\backgroundpath{
         \pgfpathmoveto{\pgfqpoint{-0.78cm}{-.17cm}}
     \pgfpathcurveto %
            {\pgfpoint{-0.35cm}{-.44cm}}
        {\pgfpoint{0.35cm}{-.44cm}}
        {\pgfpoint{.78cm}{-0.17cm}} 
     \pgfpathmoveto{\pgfqpoint{-0.78cm}{-0.17cm}}
     \pgfpathcurveto %
            {\pgfpoint{-0.25cm}{.25cm}}
        {\pgfpoint{.25cm}{.25cm}}
        {\pgfpoint{0.78cm}{-0.17cm}}
        \pgfusepath{fill}

     \pgfpathmoveto{\pgfqpoint{-1cm}{0cm}}
     \pgfpathcurveto %
            {\pgfpoint{-0.5cm}{-.5cm}}
        {\pgfpoint{0.5cm}{-.5cm}}
        {\pgfpoint{1cm}{0cm}}

     \pgfpathmoveto{\pgfqpoint{-0.75cm}{-0.15cm}}
     \pgfpathcurveto %
            {\pgfpoint{-0.25cm}{.25cm}}
        {\pgfpoint{.25cm}{.25cm}}
        {\pgfpoint{0.75cm}{-0.15cm}}
              \pgfusepath{stroke}
    }
    }

    \makeatother

\theoremstyle{plain}

\newtheorem{thm}{Theorem}[section]

\newtheorem{lem}[thm]{Lemma}

\theoremstyle{definition}
\newtheorem{dfn}[thm]{Definition}

\theoremstyle{remark}
\newtheorem{rem}[thm]{Remark}

\newtheorem{ex}[thm]{Example}

\begin{document}

\bigskip

\title{The mirror Lagrangian cobordism for the Euler exact sequence}

\date{\today}

\author{Yochay Jerby}


%
%
% ...
\begin{abstract}
For $X = \mathbb{P}^n$ the Euler sequence is given by 
$$
0 \rightarrow \Omega^1_{\mathbb{P}^n} \rightarrow \mathcal{O}_{\mathbb{P}^n}^{n+1}(-1) \rightarrow \mathcal{O}_{\mathbb{P}^n} \rightarrow 0 
$$
We describe the Lagrangian cobordism corresponding to this sequence via mirror symmetry, in the sense of Biran-Cornea. In particular, we describe the mirror Lagrangian of the cotangent sheaf $\Omega^1_{\mathbb{P}^n} \in \mathcal{D}^b(\mathbb{P}^n)$ in the mirror Fukaya category $Fuk(U_{\Delta})$. 

\end{abstract}
\maketitle
%
%
\section{Introduction and Summary of Main Results}

For $X = \mathbb{P}^n$ the Euler exact sequence is given by 
\begin{equation}
\label{eq:Eu} 
0 \rightarrow \Omega^1_{\mathbb{P}^n}  \rightarrow \mathcal{O}_{\mathbb{P}^n}^{n+1}(-1) \rightarrow \mathcal{O}_{\mathbb{P}^n} \rightarrow 0, 
\end{equation} 
see \cite{Har}. Let $\Delta $ be the moment polytope of $X$ with respect to the anti-canonical embedding. Homological mirror symmetry for toric Fano manifolds, postulates the equivalence of the categories $\mathcal{D}^b(\mathbb{P}^n)$ and $\mathcal{D}^{\pi} (Fuk(U_{\Delta}))$, where $U_{\Delta}:=Log^{-1}\vert \Delta \vert$ is the inverse image of the polytope $\Delta$ of $X$ under the logarithm map $Log \vert \cdot \vert: (\mathbb{C}^{\ast})^n \rightarrow \mathbb{R}^n$. By results of Abouzaid, mirrors of line bundles $\mathcal{O}_{\mathbb{P}^n}(k)$ are known to be given by tropical Lagrangian sections in $Fuk(U_{\Delta})$, see \cite{Ab1,Ab2}. By results of Biran-Cornea, Lagrangian cobordisims in Fukaya categories are the mirror equivalent of triangles in $\mathcal{D}^b(X)$, see \cite{BC1,BC2,BC3,BC4}. Our aim in this work is to introduce the Lagrangian cobordism corresponding to the Euler exact sequence (\ref{eq:Eu}), to which we refer as the \emph{Lagrangian Euler mirror cobordism}. Let us first consider the example of the projective line for $n=1$: 

\begin{ex}[Lagrangian Euler mirror cobordism for projective line] \label{ex:1} For $n=1$ the cotangent bundle is given by $\Omega^1_{\mathbb{P}^1}=\mathcal{O}_{\mathbb{P}^1} (-2)$. Hence, the Euler sequence is 
\begin{equation}
\label{eq:Eu1} 
0 \rightarrow \mathcal{O}_{\mathbb{P}^1} (-2) \rightarrow \mathcal{O}_{\mathbb{P}^1}(-1) \oplus \mathcal{O}_{\mathbb{P}^1}(-1) \rightarrow \mathcal{O}_{\mathbb{P}^1} \rightarrow 0. 
\end{equation} 
In this case, all three elements of the sequence are given in terms of line bundles $\mathcal{O}_{\mathbb{P}^1}(k)$. On the other hand, for the mirror, consider the annulus 
\begin{equation}
U_{[-1,1]}:=Log^{-1}([-1,1])= \left \{ e^{-1} \leq \vert z \vert \leq e \right \} \subset \mathbb{C}^{\ast}
\end{equation}
 where the map $Log \vert \cdot \vert : \mathbb{C}^{\ast} \rightarrow \mathbb{R}$ is is given by $z \mapsto Log \vert z \vert$. A Lagrangian $L$ (curve) in $U_{[-1,1]}$ is known to correspond to a line bundle $\mathcal{O}_{\mathbb{P}^1} (k)$ if it is a section of $Log \vert \cdot \vert$ beginning at $e$ and ending in $e^{-1}$. Specifically, $L$ corresponds to $\mathcal{O}_{\mathbb{P}^1} (k)$ if its lift in $\widetilde{U}_{[-1,1]}$, the universal cover of $U_{[-1,1]}$, given by $[-1,1] \times \mathbb{R}$ begins at $(e,0)$ and ends at $(e^{-1},2 \pi k)$, that is if it rotates around the origin $k \in \mathbb{Z}$ times. Applying Lagrangian surgery operations to two Lagarngians $L_1$ and $L_2$, in the sense of Polterovich \cite{P}, gives rise to a cobordism beginning in $L_1,L_2$ and ending in their surgery $L' = L_1 \# L_2$, see \cite{BC2}.
 
  Fig \ref{fig:f2} shows how the surgery of two Lagrangians $L_1$ (blue) and $L_2$ (red), which are mirror to $\mathcal{O}_{\mathbb{P}^1}(-1)\oplus \mathcal{O}_{\mathbb{P}^1}(-1)$ and $\mathcal{O}_{\mathbb{P}^1}$, results in a Lagrangian $L'=L_1 \# L_2$ mirror to $\mathcal{O}_{\mathbb{P}^1}(-2)$, giving rise to the Lagrangian Euler mirror cobordism in this case.         
\begin{figure}[ht!]
	\centering
		\includegraphics[scale=0.65]{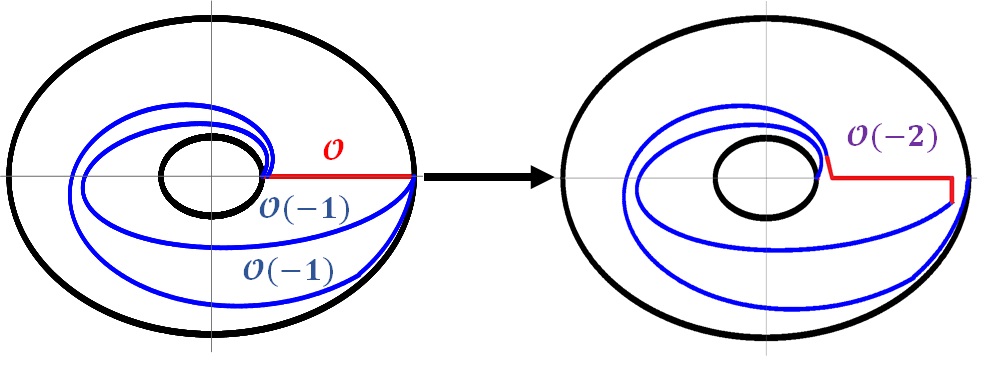} 
	\caption{Two Lagrangians corresponding to $\mathcal{O}_{\mathbb{P}^1}(-1)\oplus\mathcal{O}_{\mathbb{P}^1}(-1)$ (blue) and $\mathcal{O}_{\mathbb{P}^1}$ (red) transitioning to the Lagrangian corresponding to $\mathcal{O}_{\mathbb{P}^1}(-2)$}
\label{fig:f2}
	\end{figure}
\end{ex} 

For the higher dimensional projective spaces with $n \geq 2$ the cotangent bundle $\Omega^1_{\mathbb{P}^n}$ is no longer a line bundle but rather a vector bundle of rank $n$. In particular, the mirror of $\Omega^1_{\mathbb{P}^n}$ is not directly given by Abouzaid's construction. 
The construction of the Lagrangian Euler mirror cobordism is based on choosing two specific mirror Lagrangians $L_1,L_2$, representing $\mathcal{O}^{n+1}_{\mathbb{P}^n}(-1)$ and $\mathcal{O}_{\mathbb{P}^n}$. Our main result is Theorem \ref{thm:main} which shows that the Lagrangian $L' = L_1 \# L_2$ is the mirror of the cotangent bundle $\Omega^1_{\mathbb{P}^n}$, when considered as a $T$-equivariant bundle in the sense of Klyachko \cite{Kly}.  

\bigskip

\hspace{-0.6cm} The rest of the work is organized as follows: In section \ref{s:Cob} we recall relevant results and definitions from "both sides of the mirror". In Section \ref{s:3} we construct the Lagrangian Euler mirror cobordism for general $n \geq 1$. In Section \ref{s:4} we discuss concluding remarks.  

\section{Relevant results and definitions from both sides of the mirror} 
\label{s:Cob}

In this section we review relevant definitions and results used: 

\subsection{Projective space as a toric Fano manifold} \label{ss:2.1} A toric variety is an algebraic variety $X$ containing an algebraic torus $T \simeq (\mathbb{C}^{\ast})^n$ as a dense subset such that the action of $T$ on itself extends to the whole variety, see \cite{CLS,Fu} for standard references. A compact toric variety $X$ is said to be Fano if its anti-canonical class $-K_X$ is Cartier and ample.   

Let $N \simeq \mathbb{Z}^n$ be a lattice and let $M = N^{\vee} = Hom(N, \mathbb{Z})$ be the dual lattice. Denote by $N_{\mathbb{R}} = N \otimes \mathbb{R}$ and $M_{\mathbb{R}}= M \otimes \mathbb{R}$ the corresponding vector spaces. Let $\Delta \subset M_{\mathbb{R}}$ be an integral polytope and let 
\begin{equation} 
L(\Delta):= \bigoplus_{m \in \Delta \cap M} \mathbb{C} z^m 
\end{equation}
be the space of Laurent polynomials whose Newton polytope is $\Delta$. The polytope $\Delta$ determines an embedding 
\begin{equation}
\left \{ 
\begin{array}{c}
i_{\Delta}: (\mathbb{C}^{\ast})^n \rightarrow \mathbb{P}(L(\Delta)^{\vee}) \\ z \mapsto [z^m \vert m \in \Delta \cap M] \end{array} \right.
\end{equation}
 The polarized toric variety corresponding to the polytope $\Delta \subset \mathbb{R}^n $ is defined to be 
 \begin{equation} 
 X_{\Delta} = \overline{i_{\Delta} ((\mathbb{C}^{\ast})^n)} \subset \mathbb{P}(L(\Delta)^{\vee}),
 \end{equation} 
 the compactification of the image of $ i_{\Delta}$. The \emph{polar} polytope of $\Delta \subset M_{\mathbb{R}}$ is given by 
  \begin{equation} 
  \Delta^{\circ} = \left \{ n \vert \left < m,n \right > \geq -1 \textrm{ for every } m \in \Delta \right \} \subset N_{\mathbb{R}}  
  \end{equation} 
 The polytope $\Delta$ is said to be \emph{reflexive} if $0 \in Int (\Delta)$ and  $\Delta^{\circ}$ is an integral polytope. 
A reflexive polytope is said to be \emph{Fano} if every facet $\Delta^{\circ}$ is the convex hall of a basis of $M_{\mathbb{R}}$. Batyrev showed in \cite{Ba2} that $X_{\Delta}$ is a Fano variety if and only if $\Delta$ is reflexive and, in this case, the embedding $i_{\Delta}$ is the anti-canonical embedding. The Fano variety $X_{\Delta}$ is smooth if and only if $\Delta^{\circ}$ is a Fano polytope. The duality between the polytope $\Delta$ and its polar $\Delta^{\circ}$ serves as the basis of mirror symmetry for toric Fano manifolds (see \ref{ss:Ab-linebundles}).

Let $\Sigma=\Sigma(\Delta)$ be the fan determined by the polytope $\Delta$, see \cite{CLS,Fu}. We say that a function $\psi : N_{\mathbb{R}} \rightarrow \mathbb{R}$ is a $\Sigma$-support function if it is continuous, linear when restricted to each maximal cone $\sigma \in \Sigma(n)$, and $\psi(n_{\rho})\in \mathbb{Z}$ for any primitive generator $n_{\rho} \in N_{\mathbb{R}}$ of a one-dimensional ray $\rho \in \Sigma(1)$. Denote by $SF(\Sigma)$ the group of $\Sigma$-support functions. When $X$ is smooth one has $Div_T(X) \simeq SF( \Sigma)$ by setting 
\begin{equation}
D_{\psi} := \sum_{\rho \in \Sigma(1)} \psi(n_{\rho}) V(\rho).
\end{equation} Note that when $X$ is fano the vertices of the polar polytope $\Delta^{\circ}(0)$ are exactly the primitive generators of the one-dimensional rays of the fan $\Sigma(1)$. Denote by 
\begin{equation}
\label{eq:m}
m(\psi,\sigma) \in M = Hom(N,\mathbb{Z})
\end{equation}
 the element such that $\psi(n)=\left < n ,m(\psi,\sigma)) \right >$ for $n \in \sigma$, where $\sigma \in \Sigma(n)$ is a maximal cone.

In particular, projective space $X=\mathbb{P}^n$ is given as a toric Fano polytope $X_{\Delta}$ where $\Delta$ is the polytope whose polar is given by 
\begin{equation} 
\Delta^{\circ}= Conv \left ( \left \{-\sum_{i=1}^n e_i, e_1,...,e_n \right \} \right ) \subset N_{\mathbb{R}}.
\end{equation}    
For example, for $n=2$ the polytopes $\Delta^{\circ}$ and $\Delta$ are illustrated in Fig. \ref{fig:f3}:
\begin{figure}[ht!]
	\centering
		\includegraphics[scale=0.65]{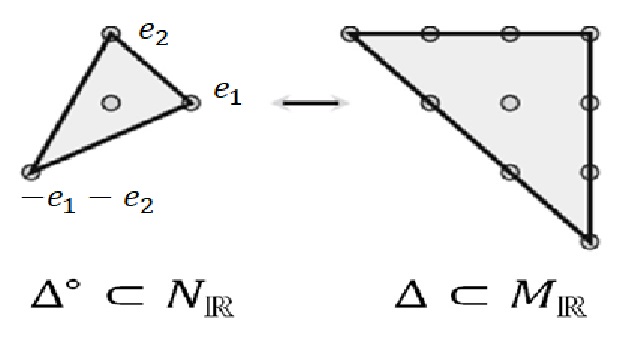} 
	\caption{The polytopes $\Delta^{\circ}$ and $\Delta$ for $X=\mathbb{P}^2$.}
\label{fig:f3}
	\end{figure}
\subsection{Line bundles and Abouzaid's tropical Lagrangian sections} \label{ss:Ab-linebundles} Let $X$ be a polarized toric manifold given by a polytope $\Delta$. Denote by $\Delta(k)$ the set of $k$-dimensional faces of $\Delta$ and let $V_X(F)$ be the closure of the $T$-orbit corresponding to the face $F \in \Delta(k)$. The group of toric divisors is given by 
\begin{equation} 
Div_T(X) := \bigoplus_{F \in \Delta(n-1)} V_X(F) \cdot \mathbb{Z}. 
\end{equation} 
The Picard group of line bundles $Pix(X)$ on a smooth toric manifold $X$ is described by the exact sequence 
\begin{equation} 
0 \rightarrow M \rightarrow Div_T(X) \rightarrow Pic(X) \rightarrow 0. 
\end{equation} 

Assume $X$ is a toric Fano manifold and let $\Delta^{\circ}$ be the polar polytope of $\Delta$.
In \cite{Ab2} Abouzaid described the mirror Lagrangian branes corresponding to elements of $Pic(X)$ as follows: Let $W \in L(\Delta^{\circ})$ be a generic Laurent polynomial whose Newton polytope is $\Delta^{\circ}$. 
Denote by $M_W=W^{-1}(0) \subset (\mathbb{C}^{\ast})^n$ the fibre of $W$ over $0 \in \mathbb{C}$. Consider 
\begin{equation} 
\widetilde{W}_{t,s}(z):=1+ \sum_{n \in \Delta^{\circ}(0)} t^{-1} (1-s \cdot \phi_n( Log \vert z \vert )z^n 
\end{equation}   
where $\phi_n \in C^{\infty}(\mathbb{R}^n)$ for $n \in \Delta^0(0)$ are required to satisfy certain decay conditions. For $s=0$ one has $\widetilde{W}_{t,0} \in L(\Delta^{\circ})$ but for $s \neq 0 $ the function $\widetilde{W}_{t,s}$ is no longer a Laurent polynomial. However, $M_{t,s}=\widetilde{W}^{-1}_{t,s}(0)$ are all symplectomorphic for generic values of $s$ for $t>>0$. Denote by $M=M_{t,1}$ with $0<<t$ big enough. We refer to the pair $((\mathbb{C}^{\ast})^n, M)$ as the \emph{tropical polarized mirror model} of the toric Fano manifold $X$. 

\begin{dfn}
A Lagrangian brane $L \subset (\mathbb{C}^{\ast})^n$ is an embedded compact graded Lagrangian submanifold, which is spin and exact. A Lagrangian brane is said to be admissible if $\partial L \subset M$ and there exists a neighbourhood of $\partial L$ in $L$ which agrees with the parallel transport of $\delta L$ along a segment $\gamma \subset \mathbb{C}$ with respect to $W$. A pair of admissible Lagrangian branes $(L_1,L_2)$ is said to be positive if their corresponding segments $\gamma_1 , \gamma_2 \subset \mathbb{C}$ lie in the half plane and their tangent vectors are oriented counter clock-wise such that $Im(\gamma_2( \theta))<Im (\gamma_1(\theta))$. 
\end{dfn}  

Admissible Lagrangian branes are objects of the Fukaya $\mathcal{A}_{\infty}$ pre-category $Fuk((\mathbb{C}^{\ast})^n, M)$ which we denote $Fuk(U_{\Delta})$, see \cite{Ab2,Ko2,S3}. The space of morphisms between two positive transverse objects $L_1,L_2 \in Fuk(U_{\Delta})$ is given by the Floer complex $(CF^{\ast}(L_1,L_2),\partial)$. In \cite{Ab2} Abouzaid introduced the $\mathcal{A}_{\infty}$ sub-pre-category of tropical Lagrangian sections 
\begin{equation}
Fuk_{trop}(U_{\Delta}) \subset Fuk(U_{\Delta})
\end{equation}  
which he proved to be quasi-equivalent to the $DG$-category of line bundles over $X$. 

Consider the map $Log \vert \cdot \vert : (\mathbb{C}^{\ast})^n \rightarrow \mathbb{R}^n$ and let $\mathcal{A} =\frac{1}{t} Log \vert M \vert  \subset \mathbb{R}^n$ be the amoeba of $M$, see \cite{GKZ}. In \cite{Ab1,Ab2} Abouzaid shows that there exists a component $\widetilde{\Delta} \subset \mathbb{R}^n \setminus \mathcal{A}$ in the complement of $\mathcal{A}$ which is contained in the polytope $\Delta \subset \mathbb{R}^n$ and is $C^0$ close to it. Note that the map $Log \vert \cdot \vert $ can be viewed as a fibration\footnote{In fact, it is considered as the dual SYZ fibration of the moment map $\mu : X \rightarrow \Delta$, in this case.} whose fibre is $\mathbb{T}^n$ and whose zero section is the Lagrangian $(\mathbb{R}^+)^n \subset (\mathbb{C}^{\ast})^n$. Consider the following definition: 

\begin{dfn} 
A tropical Lagrangian section in $((\mathbb{C}^{\ast})^n, M)$ is an admissible Lagrangian brane $L$ which is a section of the map $Log \vert \cdot \vert $ restricted to $\widetilde{\Delta}$. 
\end{dfn} 

It is shown in \cite{Ab2} that up to Hamiltonian isotopy tropical Lagrangian sections in $((\mathbb{C}^{\ast})^n, M)$ are in one-to-one correspondence with elements of $Pic(X)$. For instance, the class of the trivial bundle $\mathcal{O}_X \in Pic(X)$ corresponds to the trivial section $L_0 \subset (\mathbb{C}^{\ast})^n$, which is a tropical section of the polarized mirror model. The correspondence is based on the fact that any tropical Lagrangian section $L \subset (\mathbb{C}^{\ast})^n$ must coincide with $L_0 \subset (\mathbb{C}^{\ast})^n$ in a small neighbourhood of the fibre $Log^{-1}(x) \simeq \mathbb{T}^n$ for any vertex $x \in \Delta(0) \subset \mathbb{R}^n$. As $\mathbb{T}^n = \mathbb{R}^n / \mathbb{Z}^n$, a lift $\widetilde{L} \subset (\mathbb{R}^+)^n \times \mathbb{R}^n $ of the tropical section $L$ to the universal cover gives rise to elements $m (\widetilde{L},x) \in \mathbb{Z}^n$ for any $x \in \Delta(0)$. In particular, define $\psi_{\widetilde{L}} \in SF(\Sigma)$ to be the support function defined by 
\begin{equation} 
m(\psi, \sigma) =m(\widetilde{L},x_{\sigma}),
\end{equation}
for $x_{\sigma} \in \Delta(0) \simeq \Sigma(n)$. As a lift depends on a choice of deck transformation $m \in \mathbb{Z}^n$ we get that $ Pic(X) \simeq Fuk_{trop}(U_{\Delta}) / \sim$. 

\subsection{Klyachko's description of $T$-equivariant bundles} In \cite{Kly} Klyachko generalized the classical description of $Pic(X)$ of a toric manifold $X$, presented in \ref{ss:2.1}, to $T$-equivariant vector bundles $p: \mathcal{E} \rightarrow X$ of arbitrary rank $r = rank(\mathcal{E})$. According to Klyachko such a vector bundle is uniquely determined by a $r$-dimensional vector space $E$ equipped with filtrations $\left \{ E_{\rho}(i) \right \}_{i \in \mathbb{Z}}$ for any $\rho \in \Sigma(1)$ such that the following compatibility condition holds:  For any $\sigma \in \Sigma$ the filtrations $E_{\rho}(i)$ for $\rho \in \sigma(1)$ consist of coordinate subspaces of some basis of the space $E$. 

For instance, note that when $\mathcal{E}$ is a line bundle the vector space $E \simeq \mathbb{C}$ is a one-dimensional space. Hence, for any $\rho \in \Sigma(1)$ the filtration $E_{\rho}(i)$ is determined by the index $i_\rho \in \mathbb{Z}$ at which the filtration changes from $E$ to zero, which is the same as giving a $T$-equivariant divisor in $Div_T(X)$. 

The compatibility condition implies that the filtrations determine, for any $\sigma \in \Sigma(n)$, a decomposition $E=\bigoplus_{m \in \mathbb{Z}^n} E^{\sigma}(m)$ such that 

\begin{equation} 
E_{\rho}(i)= \sum_{\left < m, \rho \right > \geq i} E^{\sigma}(m), 
\end{equation}  

for all $\rho \in \sigma(1)$. For instance, when $\mathcal{E}= \mathcal{O}(D)$ is a line bundle with $D \in Div_T(X)$ one has 
\begin{equation} 
E^{\sigma}(m)= \left \{ \begin{array}{cc} \mathbb{C} & m=m(\psi_D, \sigma) \\ 0 & m \neq m(\psi_D, \sigma) \end{array} \right.  
\end{equation}  
Hence, the system of decompositions $E=\bigoplus_{m \in \mathbb{Z}^n} E^{\sigma}(m)$ for $\sigma \in \Sigma(n)$ generalizes the description of 
line bundles in terms of support functions. For any $T$-equivariant bundle $\mathcal{E}$ and maximal cone $\sigma \in \Sigma(n)$ let us define 
\begin{equation}
W(\mathcal{E},\sigma) := \left \{ m \vert E^{\sigma}(m) \neq 0 \right \} \subset M, 
\end{equation} 
to be the \emph{set of weights of $\mathcal{E}$ at the cone $\sigma$}. 

Klyachko gives the following description of the cotangent bundle $\Omega^1_{X}$ of a toric manifold $X$ as the vector bundle corresponding to the following system of filtrations
\begin{equation} 
\Omega_{\rho} (i) = \left \{ \begin{array}{cc} M_{\mathbb{C}} & i <0 \\ Ker \rho & i=0 \\ 0 & i >0 \end{array} \right. 
\end{equation} 
where $ker \rho = \left \{ \omega \vert \left < \omega, \rho \right >=0 \right \}$. The following example presents the corresponding decomposition determined by the filtrations $\Omega_{\rho}(i)$ for the case of projective space $X=\mathbb{P}^n$:

\begin{ex}[The weights of the cotangent bundle $\Omega^1_{\mathbb{P}^n}$]\label{ex:dec} Let $X=\mathbb{P}^n$ be realized as a toric Fano manifold by $\Delta^{\circ}(0)=\left \{ 
e_0,...,e_n \right \} $ with $e_0:=- \sum_{i=1}^n e_i$, see \ref{ss:2.1}. The maximal cones are 
\begin{equation} 
\label{eq:maximal}
\begin{array}{ccc} \sigma_0 = \sum_{i=1}^n e_i \cdot \mathbb{R}_+   & ; & \sigma_k = e_0 \cdot \mathbb{R}_+ + \sum_{i \neq k} e_i \cdot \mathbb{R}_+    \end{array} 
 \end{equation} 
 for $1 \leq k \leq n$. Hence, the decompositions of $\Omega=M_{\mathbb{C}}$, with dual basis $\left \{ e_1^{\ast},...,e_n^{\ast} \right \}$, corresponding to Klyachko's filtrations are given by 
\begin{equation} 
\label{eq:sigma10}
\Omega^{\sigma_0} = \bigoplus_{i=1}^n (-e^*_i) \cdot \mathbb{C}
  \end{equation} 
and
\begin{equation} 
\label{eq:sigma11}
\Omega^{\sigma_k} = e^*_k \cdot \mathbb{C} \oplus \left ( \bigoplus_{i \neq k} (e^*_k-e^*_i) \cdot \mathbb{C} \right ).
  \end{equation}
  In particular, the weights of $\Omega^1_X$ are given by  
  \begin{equation} 
\label{eq:sigma0}
W(\Omega^1_{\mathbb{P}^n},\sigma_0) = \left \{ -e^*_1, ...,-e^{\ast}_n \right \} 
  \end{equation} 
and
\begin{equation} 
\label{eq:sigma1}
W(\Omega^1_{\mathbb{P}^n},\sigma_k) = \left \{ e^*_k \right \} \cup \left \{ e^*_k-e^*_i \vert k \neq i \right \}.
  \end{equation}
\end{ex} 

\begin{rem}[Relation to the anti-canonical divisor] If we sum the weights of the decompositions (\ref{eq:sigma0}) and (\ref{eq:sigma1}) we get 
\begin{equation} 
m_{\sigma_0} = -e^*_1 - ... -e^*_n  
\end{equation}   
and 
\begin{equation} 
m_{\sigma_k} = n \cdot e^*_k - \sum_{i \neq k } e^*_i
\end{equation} 
for $k = 1,...,n$. Note that these are exactly the linear functionals defining the support function $\psi_{-K}$ corresponding to the anti-canonical divisor $-K_{\mathbb{P}^n}=-\sum_{i=0}^n V(\rho_i)$. That is for which $m_{\sigma_k} = m(\psi_{-K},\sigma_k)$ in the sense of (\ref{eq:m}).   

\end{rem} 

\subsection{Cobordisms, surgeries and triangles in the Fukaya category} 
\begin{figure}[ht!]
	\centering
		\includegraphics[scale=0.4]{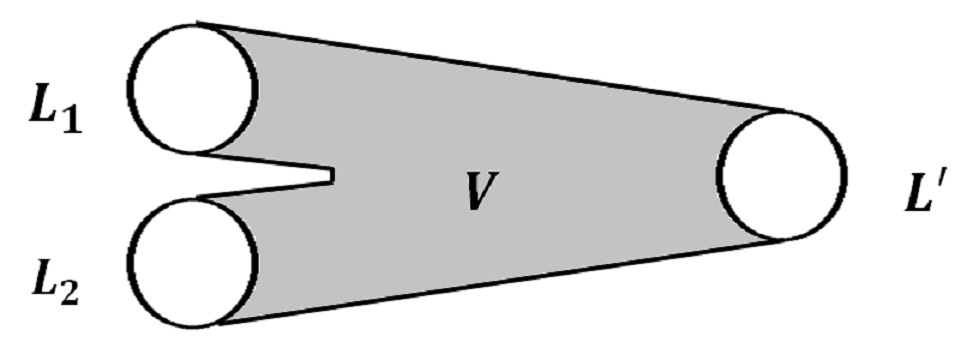} 
	\caption{A cobordism $V$ between $L_1 \cup L_2$ and $L'$.}
\label{fig:f1}
	\end{figure}

We refer the reader to the works of Biran-Cornea for the theory of Lagrangian cobordisims, see \cite{BC1,BC2,BC3,BC4}. Two Lagrangian submanifolds\footnote{ We will mainly be concerned with cobordisms with two negative ends and one positive end, as these are the cobordisms corresponding to exact sequences. In general, one can take any number of positive and negative ends.}  $L_1,L_2 \subset (M, \omega)$ are said to be Lagrangian cobordant to a third Lagrangian submanifold $L' \subset (M, \omega)$ if there exists $(V,L_1 \cup L_2, L')$, a smooth cobordism (see Fig. \ref{fig:f1}), and a Lagrangian embedding $V \subset ([0,1] \times \mathbb{R}) \times M $ such that for some $\epsilon >0$ one has 
\begin{equation}
\begin{array}{ccc}
V \mid_{[0,\epsilon) \times \mathbb{R}} = \bigcup_{i=1}^2 ([0,\epsilon) \times \left \{ i \right \} ) \times L_i) & ; & V \mid_{(1-\epsilon,1] \times \mathbb{R}} =([0,\epsilon) \times \left \{ 1 \right \} ) \times L'). 
\end{array} 
\end{equation} 

We will henceforth refer to such cobordisms as \emph{triangular}. Triangular Lagrangian cobordisms could be constructed via the method of Lagrangian surgery due Polterovich, see \cite{P}. Locally, assume $M=\mathbb{C}^n$ and consider the two Lagrangians $L_1 = \mathbb{R}^n$ and $L_2 = i \mathbb{R}^n$, intersecting transversally at the origin.  Let $H : \mathbb{R} \rightarrow \mathbb{C} $ be any smooth curve of the form $H(t)=a(t)+i b(t)$ such that 
\begin{enumerate}
\item $H(t)=t$ for $t \leq -1$
\item $H(t)=it$ for $t \geq 1$. 
\item $a'(t),b'(t)>0$ for $-1 < t < 1$. 
\end{enumerate}  
We refer to $H$ as the \emph{handle of the cobordism}. Consider the map $i_H : \mathbb{R} \times S^{n-1} \rightarrow \mathbb{C}^n$ given by 
\begin{equation} 
(t, x) \mapsto \left ( H(t) x_1,...,H(t) x_n \right ), 
\end{equation} 
where $S^{n-1}$ is considered as embedded in $\mathbb{R}^n$ with coordinates $x=(x_1,...,x_n)$. We refer to $L_H := i_H (\mathbb{R} \times S^{n-1})$ as the \emph{local surgery model} of $L_1$ and $L_2$ and denote $L_H = L_1 \#^H L_2$. We will usually omit the choice of handle $H$ and write $L' = L_1 \# L_2$. 

Globally, let $(M, \omega)$ be a general symplectic manifold and let $L_1,L_2$ be two Lagrangian sybmanifolds which intersect transversally at $L_1 \cap L_2 = \left \{ p_1,...,p_n \right \}$. One defines the surgery $L'=L_1 \# L_2$ to be a Lagrangian coinciding with $L_1 \cup L_2$ away from a small neighbourhood of the points $p_i$ and with the local surgery around $p_i$ for each $i=1,...n$, see \cite{P}.    
 
 It is known that a triangular Lagrangian cobordism $(V,L_1 \cup L_2,L')$ determines an exact triangle in $\mathcal{D}(Fuk(U_{\Delta}))$ as follows 
\begin{equation}
\begin{tikzcd}[row sep=tiny]
 L_1 \arrow[dd,"\mathcal{F}'" left] \\
& L' \arrow[ul,"\mathcal{F}" above] & \\
 L_2 \arrow[ur,"\mathcal{F}''" below] 
\end{tikzcd}
\end{equation}
where $\mathcal{F}$ is the cobordism $V$ and $\mathcal{F}',\mathcal{F}''$ are the coboridisms obtained by bending the ends of $V$ so as to turn $L'$ to the left end and $L_1$ or $L_2$, respectively, to the right end, for the general case see \cite{BC2} and \cite{FOOO}, for surgeries.  
\section{Construction of the Lagrangian Euler mirror cobordism}
\label{s:3}

Before describing the construction of the Euler mirror cobordism in the general case $n \geq 1$, let us revisit again the construction of Example \ref{ex:1} in the case $n=1$. 
\begin{ex}[The case $n=1$ via lift to the universal cover] \label{ex:2} Let $X=\mathbb{P}^1$ be the projective line and let $\Delta=[-1,1]$ be the corresponding polytope. Set 
\begin{equation} 
U_{\Delta}:= Log^{-1}(\Delta) = \left \{ e^{-1} \leq \vert z \vert \leq e \right \} \subset \mathbb{C}^{\ast}.
\end{equation}  For any $k \in \mathbb{Z}$ let $\gamma_k : I \rightarrow U_{\Delta} $ be the curve   
\begin{equation} 
\gamma_k(t) =\left ( e \cdot t + e^{-1} (1- t) \right ) e^{2 \pi k t i }.  
\end{equation} 
According to Abouzaid's mirror symmetry functor, the Lagrangian 
\begin{equation}
L(k):= \left \{ \gamma(t) \right \}_{t \in I} \subset U_{\Delta}
\end{equation} 
is the mirror representative of the line bundle $\mathcal{O}(k) \in Pic( \mathbb{P}^1)$.

Note that $ U_{\Delta} \simeq  \Delta \times \mathbb{T}$ and hence the universal cover is given by $\widetilde{U}_{\Delta} \simeq \Delta \times \mathbb{R}$. Consider the Lagrangians  
\begin{equation} 
\begin{array}{ccc} L_1 = L(-1) \cup L(-1) & ; & L_2 = L(0). \end{array} 
\end{equation} 

Corresponding to $\mathcal{O}(-1) \oplus \mathcal{O}(-1)$ and $\mathcal{O}$ of the Euler sequence. Consider the lifts $\widetilde{L}_1,\widetilde{L}_2$ of $L_1,L_2$ to the universal covering $\widetilde{U}_{\Delta} $ described in Fig. \ref{fig:f4} (the choice would be explained below, for the general case): 
\begin{figure}[ht!]
	\centering
		\includegraphics[scale=0.3]{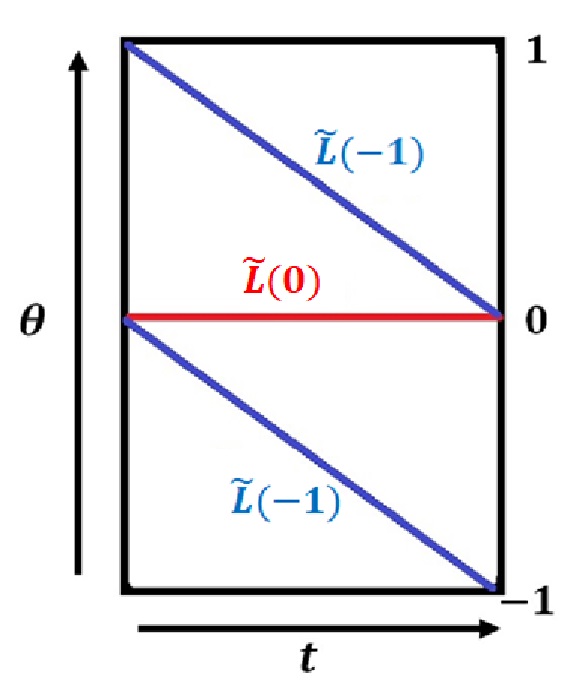} 
	\caption{The lift of $L_1 = L(1) \cup L(1)$ and $L_2 = L(0)$ to the universal cover $\widetilde{U}_{\Delta} $}
\label{fig:f4}
	\end{figure}
		
The surgery leading to $L' =L_1 \# L_2 \simeq L(-2)$ is described in Fig. \ref{fig:f5}  

	\begin{figure}[ht!]
	\centering
		\includegraphics[scale=0.55]{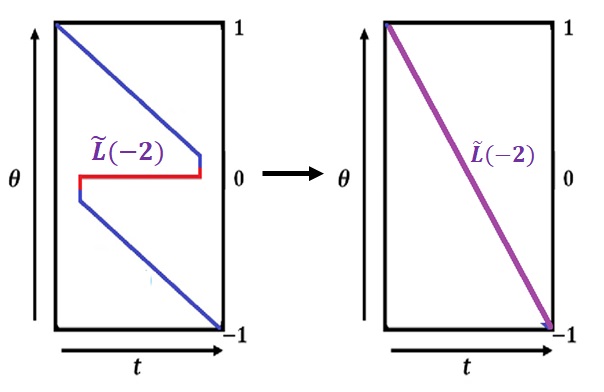} 
	\caption{The surgery $\widetilde{L}_1 \# \widetilde{L}_2$ conducted in the universal cover $\widetilde{M}$ and its equivalence to 
	$\widetilde{L}(-2)$.}
\label{fig:f5}
	\end{figure}
\end{ex}
  
In Example \ref{ex:2} we see that the lifts of the two components of $\widetilde{L}_1$ are chosen so that they each intersect the zeros section $\widetilde{L}_2$ in two different points, which are actually $(1,0) $ and $(-1,0)$, corresponding to the two vertices of the polytope $\Delta=[-1,1]$. In fact, the lift is chosen so that the components of $\widetilde{L}_1$ represent $\mathcal{O}(-D_1)$ and $\mathcal{O}(-D_{-1})$ where $Div_T (\mathbb{P}^1) = \mathbb{Z} \cdot D_1 \oplus \mathbb{Z} \cdot D_{-1}$. After applying surgery to $\widetilde{L}_1$ and $\widetilde{L}_2$ along these two points we obtain the Lagrangian $\widetilde{L}'=\widetilde{L}_1 \# \widetilde{L}_2$, which represents the line bundle whose support function $\psi$ is given by $m(\psi,\sigma_1) = -1 $ and $m (\psi, \sigma_1)=1$ (note that these are exactly the right and left heights of $\widetilde{L}'$ with respect to the $\mathbb{R}$ coordinate, as shown in Fig. \ref{fig:f5}). Direct computation shows that this is exactly the line bundle corresponding to $-K_{\mathbb{P}^1} = D_1 + D_{-1}$.     
  
 In order to generalize for any $n \geq 1$ consider the vector bundle $\mathcal{E}=\mathcal{O}^{n+1}_{\mathbb{P}^n} (-1)$ as the $T$-equivariant bundle 
 \begin{equation} 
 \mathcal{E}= \bigoplus_{i=0}^n \mathcal{O} (-V(\rho_i)),  
 \end{equation} 
 where $V(\rho_i) \in Div_T(\mathbb{P}^n)$ is the $T$-equivariant divisor corresponding to the ray $\rho_i$ for $i=0,...,n$. The following Lemma gives, via direct computation, the support functions corresponding to $V(\rho_i)$: 
 
 \begin{lem} \label{lem:1}
 Let $\psi_i := \psi_{-V(\rho_i)} \in SF(\Sigma)$ be the support-function of $-V(\rho_i)$ for $i =0,..,n$ and let $\sigma_k$ be the maximal cones of $\Sigma$ for $k=0,...,n$, as in (\ref{eq:maximal}). Then: 
\begin{enumerate}
\item For $i \neq 0$ the weights of $\psi_i$ are given by 
\begin{equation}
\begin{array}{ccccc}
m(\psi_i,\sigma_k)=e^{\ast}_k-e^{\ast}_i &  ; & m(\psi_i,\sigma_0)=-e^{\ast}_i &;  & m(\psi_i, \sigma_i)=0, 
\end{array}
\end{equation}
for $k=1,...,n$. 
\item For $i=0$ the weights of $\psi_0$ are given by 
\begin{equation} 
\begin{array}{ccc}
m(\psi_0,\sigma_k)=e^{\ast}_k &  ;  & m(\psi_0, \sigma_0)=0, 
\end{array}
\end{equation}
for $k=1,...,n$.  
\end{enumerate}
 \end{lem} 
 Recall that 
\begin{equation} 
\begin{array}{ccc} \Delta = Conv \left ( \left \{ m_{\sigma_k} \vert k=0,...,n \right \} \right ) \subset M_{\mathbb{R}} & ; & \Delta^{\circ} = Conv \left ( \left \{ e_i \vert i=0,...,n \right \} \right ) \subset  N_{\mathbb{R}}. \end{array} 
\end{equation}  
Set 
\begin{equation} 
U_{\Delta} = Log^{-1} (\Delta) \simeq \Delta \times \mathbb{T}^n \subset (\mathbb{C}^{\ast})^n.
\end{equation}
Let $\widetilde{U}_{\Delta} \simeq \Delta \times M_{\mathbb{R}}$ be the universal cover of $U_{\Delta}$ and denote by $p : \widetilde{U}_{\Delta} \rightarrow U_{\Delta}$ is the covering map. Consider the two projection maps 
\begin{center} 
\begin{tikzpicture}[commutative diagrams/every diagram]
\node (P0) at (90:2.3cm) {$\widetilde{U}_{\Delta}$};
\node (P1) at (90+72:2cm) {$\Delta$} ;
\node (P4) at (90+4*72:2cm) {$\mathbb{R}^n $};
\path[commutative diagrams/.cd, every arrow, every label]
(P0) edge node[swap] {$pr_1$} (P1)
(P0) edge node {$pr_2$} (P4);
\end{tikzpicture}
\end{center} 
Let us define: 
 
 \begin{dfn} Let $\widetilde{L} \subset \widetilde{U}_{\Delta}$ be a Lagrangian. We refer to   
 \begin{equation}
 W \left (\widetilde{L},m \right ):= pr_2 \left ( \widetilde{L}' \cap pr_1^{-1}(m) \right ) \subset \mathbb{R}^n
\end{equation}
as the \emph{set\footnote{A priori, for a general Lagrangian $\widetilde{L}$, this set is not necessarily finite, integral or non-empty.} of weights of the Lagrangian $\widetilde{L}$ at the vertex $m \in \Delta(0)$.}
 \end{dfn}

For any $i=0,...,n$ let us define 
\begin{equation} 
\widetilde{L} (-V(\rho_i)):= Conv  \left ( \left \{ (m_{\sigma_k} , m(\psi_i, \sigma_k)) \vert k=0,...,n \right \}    \right ) \subset \widetilde{U}_{\Delta},
\end{equation} 
to be the linear Lagrangian embedding of $\Delta$ in $\widetilde{U}_{\Delta}$ given by rising each vertex $m_{\sigma} \in \Delta(0)$ of $\Delta$ to height $m(\psi_i,m_{\sigma}) \in M \subset M_{\mathbb{R}}$. Finally, set 
\begin{equation} 
\begin{array}{ccc} \widetilde{L}_1 = \bigcup_{i=0}^n \widetilde{L}(-V(\rho_i)) & ; & \widetilde{L}_2 = \Delta \times \left \{ 0 \right \}, \end{array}
\end{equation} 
whose projection to $U_{\Delta}$ are mirror Lagrangians representing $\mathcal{O}_{\mathbb{P}^n}^{n+1}(-1)$ and $\mathcal{O}_{\mathbb{P}^n}$, respectively. We have: 

\begin{thm}
\label{thm:main}  
Let $\widetilde{L}'=\widetilde{L}_1 \# \widetilde{L}_2 \subset U_{\Delta} $ be the Lagrangian obtained by surgery of $\widetilde{L}_1$ and $\widetilde{L}_2$. Then 
\begin{equation} 
W \left (\widetilde{L}',m_{\sigma} \right ) = W \left (\Omega^1_{\mathbb{P}^n},\sigma \right ),
\end{equation}
for any maximal cone $\sigma \in \Sigma(n)$. In particular, $L' := p(\widetilde{L})$ is a mirror Lagrangian representing $\Omega^1_{\mathbb{P}^n}$. \end{thm}

\begin{proof} 
Note that each $\widetilde{L}(-V(\rho_i))$ intersects $\widetilde{L}_2$ uniquely in the corresponding vertex $m_{\sigma_i}$, that is 
\begin{equation} 
\widetilde{L}(\rho_i) \cap \widetilde{L}_2 = (m_{\sigma_i},0), 
\end{equation} 
for any $i = 0,...,n$. The application of surgery at this point removes from the fibre of $\widetilde{L}'$ over $m_{\sigma_i}$ the point of hight zero (compare Fig. \ref{fig:f5}). In particular, by Lemma \ref{lem:1}, for any $k=1,...,n$ one has
\begin{equation} 
W \left ( \widetilde{L}',m_{\sigma_k} \right ) = \left \{ e^{\ast}_k \right \} \cup \left \{ e_k^{\ast}-e^{\ast}_i \vert i \neq k \right \}, 
\end{equation} 
and for $k=0$ one has
\begin{equation} 
W \left ( \widetilde{L}',m_{\sigma_0} \right ) = \left \{ -e_i^{\ast} \vert i =1,...,n \right \}. 
\end{equation} 
which shows that the weights of $\widetilde{L}'$ as a Lagrangian in $\widetilde{U}_{\Delta}$ coincide with the weights of $\Omega^1_{\mathbb{P}^n}$ computed \ref{eq:sigma0} and \ref{eq:sigma1} of Example \ref{ex:dec}, as required. 
\end{proof}  

Finally, let us mention that the Euler sequence is related to mutation operations of exceptional collections, see \cite{Bo2,BoKa,GK,GR,R}:

\begin{rem}[Relation to mutation operations] \label{rem:mut} Let $\mathcal{B}= \left \{ E_1,E_2,E_3 \right \} \subset \mathcal{D}^b(\mathbb{P}^2)$ be a full strongly 
exceptional collection. The left mutation of $\mathcal{B}$ is given by 
\begin{equation} 
L(\mathcal{B})= \left \{ E_1,L_{E_2} E_3,E_2 \right \},
\end{equation} 
where the left mutation $L_{E} F$ of an exceptional object $F$ by an exceptional object $E$ is defined by the triangle 
\begin{equation} 
L_{E} F \rightarrow Hom(E, F ) \otimes E \rightarrow F \rightarrow L_{E} F[1].
\end{equation} 
According to Beilinson, the following two collections 
\begin{equation} 
\begin{array}{ccc}
\mathcal{B}_1 := \left \{ \mathcal{O}_{\mathbb{P}^2}(-1),\mathcal{O}_{\mathbb{P}^2},\mathcal{O}_{\mathbb{P}^2}(1) \right \} & ; & 
\mathcal{B}_2 := \left \{\mathcal{O}_{\mathbb{P}^2}(-1),\Omega_{\mathbb{P}^2}^1(1), \mathcal{O}_{\mathbb{P}^2} \right \} \end{array}
\end{equation} in $\mathcal{D}_{\mathbb{P}^2}^b(\mathbb{P}^n)$ are full strongly exceptional, see \cite{B}. In particular, for $E=\mathcal{O}_{\mathbb{P}^2}$ and $F=\mathcal{O}_{\mathbb{P}^2}(1)$ the left mutation $L_{\mathcal{O}_{\mathbb{P}^2}} \left (\mathcal{O}_{\mathbb{P}^2}(1) \right )$ is defined by the exact triangle  
\begin{equation} 
L_{\mathcal{O}_{\mathbb{P}^2}} \left (\mathcal{O}_{\mathbb{P}^2}(1) \right ) \rightarrow \mathcal{O}_{\mathbb{P}^2}^{n+1} \rightarrow \mathcal{O}_{\mathbb{P}^2}(1) \rightarrow L_{\mathcal{O}_{\mathbb{P}^2}} \left ( \mathcal{O}_{\mathbb{P}^2}(1) \right )[1],
\end{equation} 
which after tensoring by $\mathcal{O}_{\mathbb{P}^2}(-1)$ is exactly the Euler sequence. Hence, we have by the Euler sequence 
\begin{equation} 
L_{\mathcal{O}_{\mathbb{P}^2}} \left (\mathcal{O}_{\mathbb{P}^2}(1) \right )= \Omega^1_{\mathbb{P}^2}(1), 
\end{equation} 
and $\mathcal{B}_1= L(\mathcal{B}_2)$, that is the collection $\mathcal{B}_2$ is the left mutation of the collection $\mathcal{B}_1$. In this sense, our results can be viewed as the description of the mirror to this mutation operation in the Fukaya category $Fuk(U_{\Delta})$, we refer the reader to \cite{J,J3} where we studied the collection $\mathcal{B}_1$ from a mirror symmetry point of view.  

\end{rem}  
  
\section{Summary and concluding remarks}
\label{s:4} 

Homological mirror symmetry for toric Fano manifolds suggests the equivalence of the categories $\mathcal{D}^b(X)$ and $\mathcal{D}^{\pi} (Fuk(U_{\Delta}))$, a "dictionary" of sorts between the two categories. Currently, only a limited amount of "entries" in this dictionary are known in practice. For instance, Abouzaid's description of the mirrors of line bundles $\mathcal{O}_X(D)$ as tropical Lagrangian sections, see \cite{Ab1,Ab2}, and the description of the mirrors of structure sheaves $\mathcal{O}_{\Sigma}$ of hypersurfaces $\Sigma \subset X$ as tropical Lagrangians, see \cite{Hi,Hi2}. In this work we have described the Lagrangian cobordism mirror to the Euler short exact sequence of projective 
space $X=\mathbb{P}^n$. In particular, we obtained a description of $L'$, the mirror of the cotangent sheaf $\Omega^1_{\mathbb{P}^n}$, arising as the result of a surgery operation. This could be viewed as an extension of the "mirror dictionary" to a new type of entry, an example of the mirror of a vector bundle of rank $n$. The results actually suggest the possibility of a general framework for mirrors of (T-equivariant) vector bundles, which we hope to pursue in future work. Furthermore, the results are related to mutation operations on $\mathcal{D}^b(X)$, as explained in Remark \ref{rem:mut}.


\begin{thebibliography}{10}



\bibitem{Ab1} M.~Abouzaid. \newblock
Homogeneous coordinate rings and mirror symmetry for toric varieties \newblock
Geometry $\&$ Topology 10 (2006) 1097--1156.


\bibitem{Ab2} M.~Abouzaid. \newblock Morse homology, tropical geometry, and homological mirror symmetry for toric varieties \newblock Selecta Mathematica, August 2009, Volume 15, Issue 2, pp 189--270. 


\bibitem{Ar1} V.~I.~Arnold. \newblock Lagrange and Legendre cobordisms. \newblock I, II, Funkts. Anal. Prilozh. 14:3, 1-–13 (1980) 14:4 (1980), 8-–17


\bibitem{Ba2} V.~Batyrev. \newblock Dual polyhedra and mirror symmetry for Calabi-Yau hypersurfaces in toric varieties.
\newblock J. Algebraic Geom. 3 (1994), no. 3, 493--535.

\bibitem{BC1} P.~Biran, O.~Cornea. \newblock Lagrangian cobordism I. \newblock J. Amer. Math. Soc. 26 (2013), 295-340.


\bibitem{BC2} P.~Biran, O.~Cornea. \newblock Lagrangian cobordism and Fukaya categories. \newblock Geom. Funct. Anal. (2014), 24 (6), 1731--1830.

\bibitem{BC3} P.~Biran, O.~Cornea. \newblock Lagrangian cobordism in Lefschetz fibrations. \newblock Preprint arXiv:1504.00922 (2015)

\bibitem{BC4} P.~Biran, O.~Cornea. \newblock Cone-decompositions of Lagrangian cobordisms in Lefschetz fibrations \newblock Sel. Math. New Ser. (2017) 23: 2635--2704.


\bibitem{BC5} P.~Biran, O.~Cornea. \newblock A Lagrangian pictionary
 \newblock Kyoto J. Math. 61(2): 399-493 (June 2021). 

\bibitem{B} G. D. Birkhoff. \newblock Singular points of ordinary linear differential equations. \newblock Transactions of the American Mathematical Society. 1909, 10 (4), 436--470. 

\bibitem{B} A.~Beilinson. \newblock The derived category of coherent sheaves on $\mathbb{P}^n$. \newblock
Selected translations. Selecta Math. Soviet. 3 (1983/84), no. 3, 233--237.

\bibitem{Bo} A. ~Bondal. \newblock Helices, representations of quivers and Koszul algebras. \newblock
Helices and vector bundles, 75--95, London Math. Soc. Lecture Note Ser., 148, Cambridge Univ. Press, Cambridge, 1990.

\bibitem{Bo2} A. ~Bondal. \newblock Representations of associative algebras and coherent sheaves. \newblock
Math. USSR-Izv. 34 (1990), no. 1, 23--42.

\bibitem{BoKa} A. ~Bondal, M.~Kapranov. \newblock Representable functors, Serre functors and mutations. \newblock
Math. USSR-Izv. 35 (1990), no. 3, 519--541.

\bibitem{Bo-Or} A.~Bondal, D.~Orlov. \newblock Reconstruction of a variety from the derived category and groups of autoequivalences. \newblock Compositio Mathematica 2001, 125 (03), 327-344.

\bibitem{CLS} D. A. Cox, J. B. Little, H. K. Schenck. \newblock Toric varieties. \newblock Graduate Studies in Mathematics, 124, 2011. 


\bibitem{EP} Y.~Eliashberg, L.~Polterovich. \newblock The problem of Lagrangian knots in four-manifolds. \newblock 
Geometric Topology (Athens, 1993), AMS/IP Stud. Adv. Math., Amer. Math. Soc., 1997, 313-327. 

\bibitem{FOOO} K.~Fukaya, Y-G.~Oh, H.~Ohta, K.~Ono. \newblock Lagrangian Floer theory on compact toric manifolds. I. \newblock
Duke Math. J. 151 (2010), no. 1, 23--174.

\bibitem{FOOO3} K. Fukaya, Y-G Oh, H. Ohta, K. Ono. \newblock Lagrangian intersection Floer homology - anomaly
and obstruction - chapter 10. \newblock preprint, available at K. Fukaya's homepage (2007).

\bibitem{Fu} W. Fulton. \newblock Introduction to Toric Varieties. \newblock  Annals of Mathematics Studies, 131, Princeton University Press, 1993. 




\bibitem{GR} A. L.~Gorodentsev, A. N.~Rudakov. \newblock Exceptional vector bundles on projective spaces. \newblock Duke Math. J. 54 (1987), no. 1, 115--130. 

\bibitem{GK} A. L.~Gorodentsev, S. A.~Kuleshov. \newblock Helix theory. \newblock Mosc. Math. J., 4:2 (2004), 377--440.


\bibitem{GKZ} I. Gelfand, M. Kapranov, A. Zelevinsky. \newblock Discriminants, resultants, and multidimensional determinants. \newblock Mathematics: Theory and Applications. Birkhauser Boston, Inc., Boston, MA, 1994.

\bibitem{Har} R.~Hartshorne. \newblock Algebraic Geometry. \newblock Graduate texts in mathematics, 52, Springer-Verlag, 1977, Berlin, New York. 

\bibitem{Hi} J. Hicks. \newblock Tropical Lagrangians and Homological Mirror Symmetry. \newblock PhD Thesis, UC Berkeley

\bibitem{Hi2} J. Hicks. \newblock Tropical Lagrangian Hypersurfaces are Unobstructed. \newblock Journal of Topology 13


\bibitem{J} Y. ~Jerby. \newblock On Landau-Ginzburg systems, Quivers and Monodromy. \newblock  Journal of Geometry and Physics, 98, 2015, 504-534.

\bibitem{J3} Y. ~Jerby. \newblock On exceptional collections of line bundles and mirror symmetry for toric Del-Pezzo surfaces. \newblock Journal of Mathematical Physics 2017, 58, 031704. 

\bibitem{Kly} A.~A.~Klyachko. \newblock Equivariant Bundles on Toral Varieties. \newblock Math USSR Izv. 35 (1990), 337-375. 


\bibitem{Ko} M.~Kontsevich. \newblock Homological algebra of mirror symmetry. \newblock Proceedings of the International Congress of Mathematicians, Vol. 1, 2 (Zurich, 1994), 120-139, Birkhauser, Basel, 1995.


\bibitem{Ko2} M.~Kontsevich. \newblock Course at ENS. \newblock Preprint, 1998.

\bibitem{P} L.~Polterovich. \newblock The surgery of lagrange submanifolds. \newblock GAFA, 1991, 1 (2), 198--210. 


\bibitem{R} A. Rudakov. (Ed.) \newblock Helices and Vector Bundles: Seminaire Rudakov \newblock London Mathematical Society Lecture Note Series, 148, Cambridge University Press, 1990.


\bibitem{S3} P.~Seidel. \newblock Fukaya Categories and Picard--Lefschetz Theory, Zurich Lectures in Advanced Mathematics. \newblock European Mathematical Society (EMS), Zurich (2008)



\end{thebibliography}
\end{document}